%% file: main.tex
\documentclass[12pt]{amsart}

\usepackage{amsmath,amscd}
\usepackage{amssymb,amsthm}
\usepackage{tikz}
\usetikzlibrary{arrows.meta}
\usepackage{tikz-3dplot}

\newtheorem{theorem}{Theorem}[section]
\newtheorem{lemma}[theorem]{Lemma}

\theoremstyle{definition}

\theoremstyle{remark}
\newtheorem{remark}[theorem]{Remark}

\newcommand{\real}{\mathbb{R}}
\newcommand{\integer}{\mathbb{Z}}

\newcommand{\disk}{\mathbb{D}}
\newcommand{\bchi}{\boldsymbol{\chi}}

\begin{document}

\title[Non-stationary normal coordinates]{Non-stationary normal coordinates on neighborhoods of Pesin stable manifolds}

\author[M.~Tsujii]{Masato TSUJII}
\thanks{This work was partially supported by JSPS KAKENHI Grant Number 21H00994. The author utilized Google's Gemini for assistance with English editing, correcting typographical errors, and generating TikZ code for the figures. The author takes full responsibility for the content of the publication, of course.}
\address{Faculty of Mathematics, Kyushu University \\
Motooka 744, Nishi-ku, Fukuoka, 819-0395, Japan}
\email{tsujii@math.kyushu-u.ac.jp}

\subjclass[2020]{Primary 37D25; Secondary 37C15, 37D10}
\keywords{Non-stationary normal forms, Pesin theory, Lyapunov exponents, Invariant foliations}

\date{\today}

\begin{abstract}
We construct non-stationary normal coordinates in a neighborhood of Pesin stable manifolds \cite{Pesin}. This construction is a natural extension, via jets in the normal directions, of the non-stationary normal coordinates on stable manifolds introduced by Guysinsky and Katok \cite{GK}. We emphasize that this extension provides a useful framework for describing the fractal geometric structures of stable and unstable foliations in smooth dynamical systems.
\end{abstract}

\maketitle

\section{Introduction}
In this paper, we consider smooth dynamical systems generated by a $C^\infty$ diffeomorphism $f \colon M \to M$ on a closed Riemannian manifold $M$. Suppose that a point $p \in M$ is forward Lyapunov regular \cite{Pesin} and that there is at least one negative Lyapunov exponent at $p$. Then, Pesin's stable manifold theorem \cite{Pesin} asserts that the strong stable manifold of $p$, defined by
\[
W^s(p) = \left\{ q \in M \;\middle|\; \limsup_{n \to \infty} \frac{1}{n} \log d(f^n(q), f^n(p)) < 0 \right\},
\]
is the image of $\real^{d_s}$ under a $C^\infty$ injective immersion, where $d_s$ is the sum of the multiplicities of the negative Lyapunov exponents at $p$. That is to say, there exists a $C^\infty$ injective immersion $\iota_p \colon \real^{d_s} \to M$ such that $\iota_p(0)=p$ and $\iota_p(\real^{d_s}) = W^s(p)$.

Note that if a point $p$ is forward Lyapunov regular, then every point on its orbit is also forward Lyapunov regular. Moreover, for each $n \in \integer$, the relation $W^s(f^n(p)) = f^n(W^s(p))$ holds, and $C^\infty$ injective immersions $\iota_{f^n(p)} \colon \real^{d_s} \to M$ can be defined so that $\iota_{f^n(p)}(0) = f^n(p)$ and $\iota_{f^n(p)}(\real^{d_s}) = W^s(f^n(p))$.
The immersions $\iota_{f^n(p)}$ constructed by Pesin in \cite{Pesin} exhibit tempered distortion; namely, the distortion of $\iota_{f^n(p)}$ restricted to a compact neighborhood of the origin and that of its inverse grow sub-exponentially as $n \to \infty$. While it is natural to impose such a sub-exponential growth condition on the immersions $\iota_{f^n(p)}$, there remains considerable flexibility to modify them via composition with $C^\infty$ diffeomorphisms of $\real^{d_s}$. It is therefore a natural question whether one can choose these coordinates in an appropriate manner so as to transform the restricted action $f \colon W^s(f^n(p)) \to W^s(f^{n+1}(p))$ into a canonical normal form.
This problem is a natural extension of the classical theory of normal forms near fixed points or periodic orbits, as extensively studied by many authors including Poincar\'e and Sternberg.

The theory of non-stationary normal forms (or non-stationary normal coordinates), initiated by Guysinsky and Katok \cite{Guy,GK} and developed further by Kalinin and Sadovskaya \cite{K, KS1} and Melnick \cite{Melnick}, provides a definitive answer to this problem: we can choose the immersions $\iota_{f^n(p)} \colon \real^{d_s} \to M$ such that the representation of $f$ in these coordinates,
\[
\iota_{f^{n+1}(p)}^{-1} \circ f \circ \iota_{f^n(p)} \colon \real^{d_s} \to \real^{d_s},
\]
belongs to a finite-dimensional Lie group of polynomial automorphisms, known as sub-resonance generated polynomial automorphisms. Under mild conditions, such immersions $\iota_{f^n(p)}$ are unique up to composition with a sub-resonance generated polynomial automorphism. This implies the remarkable fact that each Pesin stable manifold carries a canonical invariant geometric structure, provided that the derivatives along the orbit satisfy certain regularity conditions. For a comprehensive review of the theory of non-stationary normal forms, we refer the reader to the survey by Kalinin \cite{Kalinin}. Regarding its various applications, there is an extensive body of literature; see the recent papers \cite{Brown, Filip, Kalinin} and the references therein.

In this paper, we discuss an extension of the aforementioned theory of non-stationary normal forms. While the existing theory has primarily focused on the dynamics \emph{within} the stable manifold itself, our motivation is to extend this framework to a neighborhood of the stable manifold so as to incorporate information concerning the normal directions. As evidenced by the seminal work of Guysinsky and Katok \cite{GK}, it is generally impossible to extend the theory of non-stationary normal forms directly to settings where both positive and negative Lyapunov exponents coexist. In this study, however, we consider the dynamics in a neighborhood of the stable manifold by treating the normal directions in terms of jets and neglecting higher-order terms. Under this formulation, the theory of non-stationary normal forms can be extended almost straightforwardly; this extension constitutes the main result of the present paper.

We emphasize that our contribution does not lie in the technical novelty of the proof itself, as our construction can be obtained by adapting the arguments of Guysinsky and Katok and subsequent studies once the setting is appropriately formulated. Rather, the significance of our work lies in the fact that the resulting normal coordinates (regarded as jets) and the associated geometric structures capture how the dynamical system ``twists'' a neighborhood of the stable manifold. This provides a useful framework for describing the intricate geometric structures of stable and unstable foliations. Indeed, the simplest case of this geometric structure has been utilized by the author in the study of exponential mixing for Anosov flows \cite{T,TZ}, and its utility has been demonstrated in several recent papers \cite{EPZ, G, GLeclerc}. We discuss further ideas for potential applications in Section \ref{s:discussion} after stating our main results.

\section{Non-autonomous dynamics on $\real^d$}
In this section, we formulate the construction of non-stationary normal coordinates for non-autonomous dynamical systems in a neighborhood of the origin in $\real^d$. As we will see in the next section, this construction immediately applies to the setting of Pesin theory and yields the main result, Theorem \ref{th:mainX}, mentioned in the introduction.

\subsection{Settings}\label{sec:setting}
We consider a sequence of $C^\infty$ diffeomorphisms defined on a bounded open neighborhood $U$ of the origin $0$ in $\real^d$:
\[
f_n \colon U \to f_n(U) \Subset \real^d, \quad n=0,1,2,\dots,
\]
satisfying $f_n(0) = 0$. We assume that the diffeomorphisms $f_n$ and their inverses
\[
f_n^{-1} \colon f_n(U) \to U, \quad n=0,1,2,\dots,
\]
are uniformly bounded for $n \ge 0$ with respect to the $C^\infty$ topology.

We assume that the linearization $L_n := (Df_n)_0$ of $f_n$ at the origin preserves a decomposition of $\mathbb{R}^d$, denoted by
\begin{equation}\label{eq:dcomE}
    \mathbb{R}^d = E_1 \oplus E_2 \oplus \cdots \oplus E_m,
\end{equation}
where $E_k \cong \mathbb{R}^{d_k}$ with $d_k \ge 1$ and $\sum_{k=1}^m d_k = d$. Furthermore, we assume that there exist a small constant $\varepsilon > 0$ and real numbers
\begin{equation}\label{eq:Chi}
\chi_1 < \chi_2 < \dots < \chi_{m_s} < 0 \le \chi_{m_s+1} < \dots < \chi_{m},
\end{equation}
such that for each $1 \le k \le m$ and $n \ge 0$, we have
\begin{equation}\label{eq:spec}
e^{\chi_k - \varepsilon} < \frac{\|L_n v\|}{\|v\|} < e^{\chi_k + \varepsilon}
\end{equation}
for any $v \in E_k \setminus \{0\}$. The constant $\varepsilon > 0$ is assumed to be sufficiently small; the precise requirements on $\varepsilon$ will be specified in Remark \ref{rem:varepsilon}.

According to the decomposition \eqref{eq:dcomE}, we sometimes denote the coordinates on $\mathbb{R}^d$ by $(x_1, \dots, x_m)$ with $x_k \in \mathbb{R}^{d_k}$. It is also natural to consider the splitting $\mathbb{R}^d = E_s \oplus E_c$ by setting
\[
E_s = E_1 \oplus \cdots \oplus E_{m_s} \quad \text{and} \quad E_c = E_{m_s+1} \oplus \cdots \oplus E_m,
\]
and to write the coordinates as $(x, y)$, where
\begin{equation}\label{eq:coordxy}
x = (x_1, \dots, x_{m_s}) \in E_s, \quad y = (x_{m_s+1}, \dots, x_m) \in E_c.
\end{equation}
In what follows, we regard the sequence $\{f_n\}_{n \ge 0}$ as a non-autonomous dynamical system, where $f_n$ maps the state at time $n$ to the state at time $n+1$.

For the arguments in the subsequent subsections, we introduce an additional assumption on $f_n$. Let $D_s$ be the closed unit disk in $E_s$ centered at the origin. We assume that the domain of definition $U$ of $f_n$ contains $D_s \times \{0\} \subset \real^d$, that
\begin{equation}\label{ass:es}
f_n(D_s \times \{0\}) \subset D_s \times \{0\},
\end{equation}
and that $(f_{n+j-1} \circ \dots \circ f_n)(D_s \times \{0\})$ contracts uniformly toward the origin as $j \to \infty$. We will consider iterations of the maps $\{f_n\}_{n\ge 0}$ as germs around the forward invariant subset $D_s \times \{0\}$.

This last assumption entails essentially no loss of generality: under the spectral conditions assumed above, the existence of a family of local stable manifolds is guaranteed by the standard theory. Thus, after rescaling if necessary, one can always perform a $C^\infty$ change of coordinates in a neighborhood of the origin so that \eqref{ass:es} and the uniform contraction on $D_s \times \{0\}$ are satisfied.

\begin{remark}
The setting described above is more restrictive than those typically found in the literature. Specifically, we assume $C^\infty$ regularity, and our spectral condition \eqref{eq:spec} is more stringent than the ``narrow band condition'' introduced in \cite{GK}. The assumption of $C^\infty$ regularity is made to simplify the exposition; the arguments presented here can be adapted to $C^r$ diffeomorphisms provided that $r$ is sufficiently large (depending on the data $\{\chi_k\}$ and $\varepsilon$). Such regularity requirements were extensively studied by Kalinin \cite{Kalinin} in the original setting of non-stationary normal forms. We adopt the stronger condition \eqref{eq:spec} to keep the presentation concise and to focus on the most relevant applications.
\end{remark}

\subsection{Coordinates and Jets}
By a $C^\infty$ closed disk $V$ in $E_s$, we mean the image $V = \varphi(D_s)$ under a $C^\infty$ diffeomorphism $\varphi \colon E_s \to E_s$ with $\varphi(0)=0$. For a $C^\infty$ closed disk $V$ in $E_s$ or for $V=E_s$, let $\mathcal{C}^{\infty}(V)$ be the set of $C^\infty$ maps $h$ from a neighborhood of $V \times \{0\}$ in $\real^d = E_s \oplus E_c$ to $\real^d$ such that $h(0)=0$ and $h(V \times \{0\}) \subset E_s \times \{0\}$.
Let $\mathcal{C}_0^\infty(V)$ be the subspace consisting of maps $h \in \mathcal{C}^{\infty}(V)$ that additionally satisfy $(Dh)_0 = 0$.
For each $h \in \mathcal{C}^\infty_{0}(V)$, the map $H = \mathrm{Id} + h$ is a $C^\infty$ diffeomorphism when restricted to a sufficiently small neighborhood of the origin. We will therefore regard such a map $H$ as a coordinate change defined near the origin.

Recall the coordinates $(x,y)$ introduced in \eqref{eq:coordxy}. Let $\ell \ge 0$ be an integer. Two elements $g_1, g_2 \in \mathcal{C}^\infty(V)$ are said to be \textbf{$\ell$-jet equivalent in the directions transversal to $E_s$}, denoted by $g_1 \sim_\ell g_2$, if
\[
\|g_1(x,y) - g_2(x,y)\| = o(\|y\|^\ell) \quad \text{as } y \to 0
\]
for each $x \in V$. Equivalently, $g_1 \sim_\ell g_2$ for $g_1,g_2 \in \mathcal{C}^\infty(V)$ if and only if
\[
\partial_y^{\beta} g_1(x,0) = \partial_y^{\beta} g_2(x,0) \quad \text{for all } x \in V \text{ and multi-indices } |\beta| \le \ell.
\]
We denote by $[g]^{(\ell)}$ the equivalence class of $g \in \mathcal{C}^\infty(V)$ and consider the quotient spaces
\begin{equation}\label{def:Cell}
\mathcal{C}^{(\ell)}(V) := \mathcal{C}^\infty(V) / {\sim_{\ell}} \quad \text{and} \quad \mathcal{C}^{(\ell)}_0(V) := \mathcal{C}^\infty_0(V) / {\sim_{\ell}}.
\end{equation}

\subsection{Sub-resonance generated polynomial automorphisms}
\label{ss:subr}

For an integer $\ell \ge 0$, let $\mathcal{P}^{(\ell)}_0$ be the space of polynomial maps $p \colon \real^d \to \real^d$ whose degree with respect to the variables $y \in E_c$ is at most $\ell$, and which satisfy
\[
p(0)=0, \quad (Dp)_0 = 0, \quad \text{and} \quad p(E_s \times \{0\}) \subset E_s \times \{0\}.
\]
For a $C^\infty$ closed disk $V$ in $E_s$ or for $V=E_s$, there is a natural injection
\begin{equation}\label{eq:id}
[\cdot ]^{(\ell)} \colon \mathcal{P}^{(\ell)}_0 \to \mathcal{C}^{(\ell)}_0(V), \quad p \mapsto [p]^{(\ell)}.
\end{equation}
We will identify $\mathcal{P}^{(\ell)}_0$ and its subsets with their images in $\mathcal{C}^{(\ell)}_0(V)$ via this injection.

\begin{figure}[ht]
    \centering
    \input{pict1.tex}
    \caption{The action of $\mathcal{L}_{\bchi}$ around $0$.}
    \label{fig1}
\end{figure}
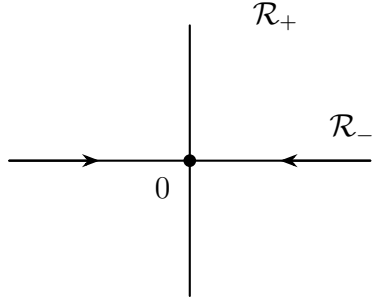

The following definitions refer to the decomposition \eqref{eq:dcomE} and the spectral data \eqref{eq:Chi}, which we collectively denote by $\bchi$.
Let $L_{\bchi} \colon \real^d \to \real^d$ be the unique linear isomorphism specified by
\begin{equation}\label{eq:Lchi}
    L_{\bchi}(v) = e^{\chi_k} v \quad \text{for all } v \in E_k, \;\; 1 \le k \le m.
\end{equation}
We define an action of $L_{\bchi}$ on $\mathcal{P}^{(\ell)}_0$ via conjugation:
\[
\mathcal{L}_{\bchi} \colon \mathcal{P}_0^{(\ell)} \to \mathcal{P}^{(\ell)}_0, \quad p \mapsto L_{\bchi}^{-1} \circ p \circ L_{\bchi}.
\]
The operator $\mathcal{L}_{\bchi}$ is diagonalizable with real eigenvalues. Consequently, the space $\mathcal{P}_0^{(\ell)}$ admits a direct sum decomposition
\begin{equation}\label{eq:Pdecomp}
\mathcal{P}^{(\ell)}_0 = \mathcal{R}^{(\ell)}_+ \oplus \mathcal{R}^{(\ell)}_-,
\end{equation}
where $\mathcal{R}^{(\ell)}_+$ (resp.\ $\mathcal{R}^{(\ell)}_-$) is the finite-dimensional (resp.\ infinite-dimensional) subspace spanned by the eigenvectors of $\mathcal{L}_{\bchi}$ corresponding to eigenvalues $\lambda \ge 1$ (resp.\ $\lambda < 1$). See Figure \ref{fig1}.

In Section \ref{ss:pf}, we will provide a more precise description of this decomposition and prove the following lemma.

\begin{lemma}\label{lm:SR}
The set
\[
\mathcal{G}^{(\ell)}_+ := \{ [\mathrm{Id} + p]^{(\ell)} \mid p \in \mathcal{R}^{(\ell)}_+ \} \subset \mathcal{C}^{(\ell)}(E_s)
\]
forms a finite-dimensional Lie group with respect to composition. Similarly, the set
\[
\mathcal{G}^{(\ell)}_- := \{ [\mathrm{Id} + p]^{(\ell)} \mid p \in \mathcal{R}_-^{(\ell)} \} \subset \mathcal{C}^{(\ell)}(E_s)
\]
is closed under composition. (Note that the inverse of an element of $\mathcal{G}^{(\ell)}_-$ cannot be represented by a polynomial map in general.)
\end{lemma}

Let $\mathcal{Z}$ be the group of linear automorphisms of $\real^d$ that preserve the decomposition \eqref{eq:dcomE}.
Let $\mathcal{Z}'$ be the group of linear automorphisms of $\real^d$ that preserve the following flag of linear subspaces of $\real^d$:
\begin{equation}\label{eq:flag}
E_1 \subset E_1 \oplus E_2 \subset E_1 \oplus E_2 \oplus E_3 \subset \cdots \subset \bigoplus_{k=1}^{m} E_k = \real^d.
\end{equation}
We then consider the semi-direct products
\begin{align*}
 \mathcal{G}^{(\ell)}_+ \rtimes \mathcal{Z} &= \{ [(\mathrm{Id} + p) \circ B]^{(\ell)} \mid p \in \mathcal{R}^{(\ell)}_+, \; B \in \mathcal{Z} \} \subset \mathcal{C}^{(\ell)}(E_s), \quad \text{and}\\
 \mathcal{G}^{(\ell)}_+ \rtimes \mathcal{Z}' &= \{ [(\mathrm{Id} + p) \circ B]^{(\ell)} \mid p \in \mathcal{R}^{(\ell)}_+, \; B \in \mathcal{Z}' \} \subset \mathcal{C}^{(\ell)}(E_s).
\end{align*}
It is easy to check that these indeed form groups with respect to composition.
We call the elements of these groups \textbf{sub-resonance generated polynomial automorphisms}.

\subsection{Main Result}
We now formulate our main result for the non-autonomous dynamical system $\{f_n \colon U \to \real^d\}_{n \ge 0}$ satisfying the assumptions given in Section \ref{sec:setting}. We defer the proof to the last section.

\begin{theorem}\label{th:abs}
There exist an open neighborhood $W$ of $D_s \times \{0\}$ in $\real^d$ and a sequence of maps $\{h_n \colon W \to \real^d\}_{n \ge 0}$ in $\mathcal{C}^{\infty}_0(D_s)$, which is uniformly bounded in the $C^\infty$ topology, such that for all $n \ge 0$, the $\ell$-jet equivalence classes of the coordinate representations of $f_n$,
\[
P_n := [(\mathrm{Id} + h_{n+1}) \circ f_n \circ (\mathrm{Id} + h_n)^{-1}]^{(\ell)} \in \mathcal{C}^{(\ell)}(V_n),
\]
belong to the Lie group $\mathcal{G}^{(\ell)}_+ \rtimes \mathcal{Z}$, where $V_n = (\mathrm{Id}+h_n)(D_s) \subset E_s$.
Furthermore, among sequences $\{h_n\}_{n \ge 0}$ that are uniformly bounded in the $C^\infty$ topology, the equivalence class $[\mathrm{Id} + h_n]^{(\ell)} \in \mathcal{C}^{(\ell)}(D_s)$ for each $n \ge 0$ is unique up to composition on the left with elements of $\mathcal{G}^{(\ell)}_+$.
\end{theorem}

In the case $\ell=0$, the theorem above provides a sequence of coordinate changes on $D_s$, which essentially recovers the conclusions of the standard theory of non-stationary normal forms \cite{GK, KS1, Melnick}. In the case $\ell > 0$, our result yields additional geometric information regarding the behavior of the dynamics in the directions transversal to the stable manifold, as noted in the introduction.

\begin{remark}
Since we consider non-autonomous dynamics in forward time $n \ge 0$, the asymptotic behavior of the dynamics itself determines not the individual subspaces $E_k$ ($1 \le k \le m$), but rather the flag \eqref{eq:flag}. If one wishes to take this ambiguity into account, one can start with the flag \eqref{eq:flag} instead of the decomposition \eqref{eq:dcomE} and obtain the same conclusion with the group $\mathcal{G}^{(\ell)}_+ \rtimes \mathcal{Z}'$ in place of $\mathcal{G}^{(\ell)}_+ \rtimes \mathcal{Z}$.
\end{remark}

\section{An application to the framework of Pesin theory}
We now apply the results of the previous section to the framework of Pesin theory to obtain the main theorem mentioned in the introduction.

Let us recall a few basic notions from smooth ergodic theory. Consider a dynamical system generated by a $C^\infty$ diffeomorphism $f \colon M \to M$ on a closed manifold $M$. We fix a Riemannian metric on $M$. For each $p \in M$, the forward Lyapunov exponent is defined as a function on the tangent space $T_pM \setminus \{0\}$ by
\[
\chi(v) = \limsup_{n \to \infty} \frac{1}{n} \log \|Df_p^n(v)\| \quad \text{for } v \in T_pM \setminus \{0\}.
\]
This function takes finitely many values at each point $p$:
\begin{equation}\label{eq:le}
\chi_1(p) < \chi_2(p) < \dots < \chi_{m(p)}(p),
\end{equation}
with corresponding multiplicities $d_1(p), d_2(p), \dots, d_{m(p)}(p) \ge 1$ satisfying $\sum_{k=1}^{m(p)} d_k(p) = \dim M$. These multiplicities are given by the dimensions of the subspaces in the filtration:
\[
\dim\{ v \in T_pM \mid \chi(v) \le \lambda \} = \sum_{k\colon \chi_k(p) \le \lambda} d_k(p) \quad \text{for any } \lambda \in \real.
\]
To avoid trivial cases, we assume that there is at least one negative Lyapunov exponent at $p$.

From the sub-additivity of the growth rates, we always have
\begin{equation}\label{eq:Lya}
\sum_{k=1}^{m(p)} \chi_k(p) d_k(p) \ge \liminf_{n \to \infty} \frac{1}{n} \log |\det Df_p^n|.
\end{equation}
A point $p \in M$ is said to be \textbf{forward Lyapunov regular} if equality holds in \eqref{eq:Lya}. By the Oseledets Multiplicative Ergodic Theorem \cite{Oseledec}, almost every point with respect to any $f$-invariant Borel probability measure is forward Lyapunov regular.

For a forward Lyapunov regular point $p \in M$, there exists a splitting (the Oseledets decomposition) of the tangent space:
\[
T_pM = E_1(p) \oplus E_2(p) \oplus \dots \oplus E_{m(p)}(p),
\]
where $\dim E_k(p) = d_k(p)$, such that for any $0 \neq v \in E_k(p)$, the limit
\[
\lim_{n \to \infty} \frac{1}{n} \log \|Df_p^n(v)\| = \chi_k(p)
\]
exists. Moreover, for any $k \neq k'$, the angle between the subspaces satisfies
\[
\lim_{n \to \infty} \frac{1}{n} \log \angle(Df^n_p(E_k(p)), Df^n_p(E_{k'}(p))) = 0.
\]
Adapting the argument from the previous section, we obtain the following theorem.

\begin{theorem}\label{th:mainX}
Suppose that $p \in M$ is a forward Lyapunov regular point, and let the Lyapunov exponents $\{\chi_k = \chi_k(p)\}_{k=1}^m$ at $p$ and their multiplicities $\{d_k = d_k(p)\}_{k=1}^m$ be given as above.
Let $1 \le m_s \le m$ be such that \eqref{eq:Chi} holds, and assume that $\real^d$ is decomposed as in \eqref{eq:dcomE}.
Then there exists a sequence of $C^\infty$ embeddings
\[
\hat{\iota}_{f^n(p)} \colon U \to M
\]
defined on an open neighborhood $U$ of $D_s \times \{0\}$, satisfying $\hat{\iota}_{f^n(p)}(0) = f^n(p)$, $\hat{\iota}_{f^n(p)}(D_s \times \{0\}) \subset W^s(f^n(p))$, and $(D\hat{\iota}_{f^n(p)})_0(E_k) = Df^n_p(E_k(p))$, such that the following hold:
\begin{itemize}
    \item[(i)] The $\ell$-jets of the coordinate representations of $f$, given by
    \[
    \left[\hat{\iota}_{f^{n+1}(p)}^{-1} \circ f \circ \hat{\iota}_{f^n(p)}\right]^{(\ell)} \in \mathcal{C}^{(\ell)}(D_s) \qquad \text{for } n \ge 0,
    \]
    belong to the Lie group $\mathcal{G}^{(\ell)}_+ \rtimes \mathcal{Z}$.
    \item[(ii)] The maps $\hat{\iota}_{f^n(p)}$ (restricted to $U$) and their inverses grow sub-exponentially as $n \to \infty$ in the $C^\infty$ topology.
\end{itemize}
Furthermore, such a sequence $\{\hat{\iota}_{f^n(p)}\}_{n \ge 0}$ is unique up to $\ell$-jet equivalence and composition with elements of $\mathcal{G}^{(\ell)}_+ \rtimes \mathcal{Z}$.
\end{theorem}

\begin{remark}
Condition (ii) above is stated precisely as follows.
For $n \ge 0$, let $\disk(\gamma(n)) \subset \real^d$ be the open disk of radius $\gamma(n) > 0$ centered at the origin. Provided that $\gamma(n)$ is sufficiently small, the image of $\disk(\gamma(n))$ under $\hat{\iota}_{f^n(p)}$ is contained in a small neighborhood of $f^n(p)$, allowing us to measure the $C^r$ norms of the restriction $\hat{\iota}_{f^n(p)}|_{\disk(\gamma(n))}$ and its inverse via local coordinates on $M$. Condition (ii) requires the existence of a sequence $\gamma(n) > 0$ decaying at most sub-exponentially as $n \to \infty$ such that, for every $r \ge 1$, the $C^r$ norms of $\hat{\iota}_{f^n(p)}|_{\disk(\gamma(n))}$ and its inverse grow sub-exponentially as $n \to \infty$.
\end{remark}

\section{Discussion}\label{s:discussion}
Before proceeding to the proofs of the results stated in the previous sections, we outline potential applications of the non-stationary normal coordinates constructed in Theorem~\ref{th:mainX} to the geometric analysis of stable and unstable foliations in smooth dynamical systems.

\begin{figure}[ht]
    \centering
    \input{pict2.tex}
    \caption{Local stable manifolds for points on $W^u_{\mathrm{loc}}(p)$.}
    \label{fig:localstable}
\end{figure}
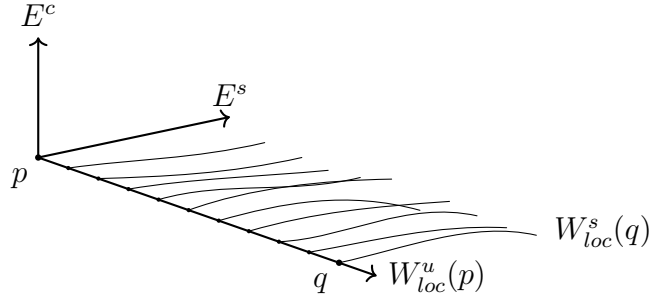

Consider a dynamical system generated by a sufficiently smooth diffeomorphism $f \colon M \to M$ on a closed manifold $M$. For simplicity, assume that the volume measure is ergodic and that $f$ possesses both positive and negative Lyapunov exponents. Under these assumptions, for almost every point $p \in M$, there exists a local unstable manifold $W^u_{\mathrm{loc}}(p)$. We are interested in the family of local stable manifolds $\{W^s_{\mathrm{loc}}(q)\}_{q \in W^u_{\mathrm{loc}}(p)}$ at points along the unstable manifold of $p$ (see Figure \ref{fig:localstable}).

In general, the dependence of $W^s_{\mathrm{loc}}(q)$ on the base point $q \in W^u_{\mathrm{loc}}(p)$ is known to be only H\"older continuous rather than smooth. This lack of regularity is a fundamental manifestation of the non-smoothness of the holonomy maps associated with the stable foliation. When a foliation fails to be smooth, classical differential-geometric arguments are difficult to apply directly. Nevertheless, it remains possible to study the structure of the stable foliation from the viewpoint of fractal geometry or dynamical renormalization.

Indeed, if we apply the iterates $f^n$ to the family $\{W^s_{\mathrm{loc}}(q)\}_{q \in W^u_{\mathrm{loc}}(p)}$, we obtain a new family of local stable manifolds $\{f^n(W^s_{\mathrm{loc}}(q))\}$ associated with the points $f^n(q)$, which lie on the local unstable manifold $W^u_{\mathrm{loc}}(f^n(p))$. In this sense, the family of local stable manifolds exhibits a form of non-linear self-similarity under the action of the dynamics.

This observation becomes significantly more powerful when combined with the results of the present paper. By applying Theorem~\ref{th:mainX} to the inverse map $f^{-1}$, we can choose appropriate coordinates along the unstable manifolds such that the transformation of a neighborhood---at the level of $\ell$-jets---is described by elements of a finite-dimensional Lie group of sub-resonance generated polynomial automorphisms.

We anticipate that this framework, which captures the infinitesimal ``twisting'' of neighborhoods via the language of jets, will provide a robust tool for the analysis of stable and unstable foliations. In particular, we hope that these normal coordinates will facilitate a deeper understanding of the geometric rigidity and statistical properties of smooth dynamical systems with uniform or non-uniform hyperbolicity.

Finally, we note that the ideas outlined above have already been applied in the author's earlier papers \cite{T,TZ} in the more restrictive setting of partially hyperbolic dynamical systems in dimension three. In that context, we used normal coordinates on the conormal bundle of the stable manifolds, which are naturally induced by the normal coordinates on the jet space. The functions describing the fluctuations of the unstable direction along stable manifolds are called \emph{template functions}; collections of such functions are useful for quantifying the joint non-integrability of the stable and unstable foliations. See \cite{T,TZ} for details.

\section{Proofs}
\subsection{Proof of Lemma \ref{lm:SR}} \label{ss:pf}
We first examine the subspace $\mathcal{R}^{(\ell)}_+$ in detail. An element of $\mathcal{P}^{(\ell)}_0$ is represented by a polynomial map whose $k$-th component $p_k \colon \real^d \to E_k$ (for $1 \le k \le m$) is written as
\begin{equation}\label{eq:pk}
p_k(x,y) = \sum_{(\alpha,\beta)} C_{k,\alpha,\beta} \; x^{\alpha} y^{\beta}, \qquad C_{k,\alpha,\beta} \in E_k,
\end{equation}
where the sum $\sum_{(\alpha,\beta)}$ is taken over multi-indices
\[
(\alpha,\beta) = ((\alpha_1, \dots, \alpha_{m_s}), (\alpha_{m_s+1}, \dots, \alpha_{m})) \in \integer_{\ge 0}^{d_s} \times \integer_{\ge 0}^{d-d_s}
\]
with $\alpha_j \in \integer_{\ge 0}^{d_j}$ ($1 \le j \le m$) satisfying $|\beta| \le \ell$ and $|\alpha| + |\beta| \ge 2$ (and $|\beta| \ge 1$ if $k > m_s$), where $d_s = \sum_{k=1}^{m_s} d_k$. Each monomial term $C_{k,\alpha,\beta} x^{\alpha} y^{\beta}$ in the expression for $p_k$ is an eigenvector of the operator $\mathcal{L}_{\bchi}$ with the corresponding eigenvalue
\[
\exp\left( -\chi_k + \sum_{j=1}^m |\alpha_j| \chi_j \right).
\]
Since $\chi_j < 0$ for $j \le m_s$, this eigenvalue is strictly less than $1$ for all but finitely many pairs of $1 \le k \le m$ and multi-indices $(\alpha,\beta)$ with $|\beta| \le \ell$.
Consequently, the subspace $\mathcal{R}^{(\ell)}_+$ is finite-dimensional.
Moreover, one can verify by direct computation that the composition of two elements of $\mathcal{G}^{(\ell)}_+$ remains in $\mathcal{G}^{(\ell)}_+$. Furthermore, by considering the inversion process at the level of formal power series (modulo $o(\|y\|^\ell)$), it follows that the inverse of any element in $\mathcal{G}^{(\ell)}_+$ also belongs to $\mathcal{G}^{(\ell)}_+$ (cf.\ \cite[Lemma 1.1]{GK}). This completes the proof of Lemma \ref{lm:SR} for $\mathcal{G}^{(\ell)}_+$. The claim for $\mathcal{G}^{(\ell)}_-$ follows from a similar argument.

\subsection{Proof of Theorem \ref{th:abs}} \label{ss:pf_main}
We begin with a preparatory reduction.
Let $f_n \colon U \to \real^d$ be the sequence of maps in the statement of the theorem.
Consider the rescaled maps
\[
\tilde{f}_n = \Lambda \circ f_n \circ \Lambda^{-1} \colon \Lambda(U) \to \real^d \quad \text{for } n \ge 0,
\]
where $\Lambda$ denotes scalar multiplication by a large factor $\lambda \ge 1$.
Suppose that the conclusion of the theorem holds for the rescaled sequence $\{\tilde{f}_n\}$ with the corresponding maps $\tilde{h}_n$ and $\tilde{P}_n$.
Then we obtain the conclusion for the original sequence $\{f_n\}$ by setting $P_n = \Lambda^{-1} \circ \tilde{P}_n \circ \Lambda$ and defining $h_n$ as
\[
h_n = (P_{n+j} \circ \dots \circ P_n)^{-1} \circ \Lambda^{-1} \circ (\mathrm{Id} + \tilde{h}_{n+j+1}) \circ \Lambda \circ (f_{n+j} \circ \dots \circ f_n) - \mathrm{Id}
\]
for a sufficiently large integer $j \ge 0$.
This reduction implies that, in proving the theorem, we may assume without loss of generality that the non-linear part $f_n - (Df_n)_0$ is arbitrarily small in the $C^\infty$ topology, uniformly in $n \ge 0$.

Let $r \ge 2$ be a sufficiently large integer such that the degree with respect to the variables $x$ of any component of any element in $\mathcal{R}^{(\ell)}_+$ is at most $r$, and let $R \ge r$ be an arbitrarily large integer.
Let $V$ be a $C^\infty$ closed disk in $E_s$. Each $\ell$-jet equivalence class $[h]^{(\ell)} \in \mathcal{C}^{(\ell)}_0(V)$ has a unique canonical representative of the form
\begin{equation}\label{eq:hxy}
h(x,y) = \sum_{\beta \colon |\beta|\le \ell} y^{\beta} \cdot h_{\beta}(x) \quad \text{in } \mathcal{C}^\infty_0(V).
\end{equation}
Using this representation, we define the $C^R$ norm $\|\cdot\|_{C^R}$ on $\mathcal{C}^{(\ell)}_0(V)$ by
\[
\|[h]^{(\ell)}\|_{C^R} = \max_{|\beta|\le \ell} \; \max_{|\alpha|\le R} \;
\sup_{x\in V} \|\partial^\alpha_x h_{\beta}(x)\|.
\]
Now consider a linear transformation $L \colon \real^d \to \real^d$ that preserves the decomposition \eqref{eq:dcomE} and satisfies \eqref{eq:spec} for a small $\varepsilon > 0$.
Its action on the space $\mathcal{C}^{(\ell)}_0(L(V))$ by conjugation is denoted by
\begin{equation}\label{eq:L}
\mathcal{L} \colon \mathcal{C}^{(\ell)}_0(L(V)) \to \mathcal{C}^{(\ell)}_0(V), \quad \mathcal{L}([h]^{(\ell)}) = [L^{-1} \circ h \circ L]^{(\ell)}.
\end{equation}
Note that this action is analogous to that of $\mathcal{L}_{\bchi}$ on the space $\mathcal{P}^{(\ell)}_0$ discussed in Sections \ref{ss:subr} and \ref{ss:pf}. Indeed, $\mathcal{L}$ preserves the direct sum decomposition analogous to \eqref{eq:Pdecomp}:
\begin{equation}\label{eq:prod}
\mathcal{C}^{(\ell)}_0(V) = \mathcal{R}^{(\ell)}_+ \oplus \overline{\mathcal{R}}_{-}^{(\ell)}(V)
\end{equation}
and the corresponding decomposition of $\mathcal{C}^{(\ell)}_0(L(V))$, where
\[
\overline{\mathcal{R}}_{-}^{(\ell)}(V) = \left\{ p + [h]^{(\ell)} \;\middle|\;
p \in \mathcal{R}^{(\ell)}_- \text{ and } \partial^{\alpha}_x h_{\beta}(0) = 0 \text{ for all } |\alpha| \le r, \; |\beta| \le \ell \right\},
\]
in which $h \in \mathcal{C}_0^\infty(V)$ is taken to be of the form \eqref{eq:hxy}.
Furthermore, the action of $\mathcal{L}$ on the second component,
\begin{equation}\label{eq:second}
\mathcal{L} \colon \overline{\mathcal{R}}_{-}^{(\ell)}(L(V)) \to \overline{\mathcal{R}}_{-}^{(\ell)}(V),
\end{equation}
is a contraction. More precisely, we can modify the $C^R$ norm $\|\cdot\|_{C^R}$ on $\mathcal{C}^{(\ell)}_0(V)$ within its equivalence class of norms so that, for some constant $0 < \rho_R < 1$, we have
\begin{equation}\label{eq:contraction}
\|\mathcal{L} \colon \overline{\mathcal{R}}_{-}^{(\ell)}(L(V)) \to \overline{\mathcal{R}}_{-}^{(\ell)}(V)\|_{C^R} \le \rho_R < 1
\end{equation}
for any linear transformation $L$ that preserves the decomposition \eqref{eq:dcomE} and satisfies \eqref{eq:spec} for a sufficiently small $\varepsilon > 0$.
In addition, we may and do assume that the projection $\pi \colon \mathcal{C}^{(\ell)}_0(V) \to \overline{\mathcal{R}}_{-}^{(\ell)}(V)$ along the subspace $\mathcal{R}^{(\ell)}_+$ has operator norm $1$ with respect to this modified $C^R$ norm $\|\cdot\|_{C^R}$.

\begin{remark}\label{rem:varepsilon}
For the argument above to hold, we must choose the constant $\varepsilon > 0$ in assumption \eqref{eq:spec} to be sufficiently small. Specifically, the requirement is that whenever $\chi'_j \in [\chi_j - \varepsilon, \chi_j + \varepsilon]$ for $1 \le j \le m$, we have
\[
-\chi_k + \sum_{j=1}^m |\alpha_j| \chi_j < 0 \quad
\Longleftrightarrow \quad
-\chi'_k + \sum_{j=1}^m |\alpha_j| \chi'_j < 0
\]
for any $1 \le k \le m$ and any $(\alpha_j)_{j=1}^m = (\alpha,\beta) \in \mathbb{Z}_{\ge 0}^{d_s} \times \mathbb{Z}_{\ge 0}^{d-d_s}$ with $|\beta| \le \ell$ and $|\alpha| + |\beta| \ge 2$.
\end{remark}

We now proceed to the proof of Theorem \ref{th:abs}. We represent the maps $f_n \colon U \to \real^d$ in the form
\[
f_n = (\mathrm{Id} + g_n) \circ L_n \quad \text{with } g_n \in \mathcal{C}^\infty_0(L_n(D_s)),
\]
where $L_n = (Df_n)_0$ is the linearization of $f_n$ at the origin. By the preliminary reduction at the beginning of this subsection, we may and do assume that
\begin{equation}\label{eq:gn}
\|[g_n]^{(\ell)}\|_{C^R} < \delta \quad \text{uniformly in } n \ge 0
\end{equation}
for a sufficiently small constant $\delta > 0$.
For each integer $N \ge 2$, we define two finite sequences of maps:
\begin{itemize}
    \item $H_n^N = \mathrm{Id} + h_n^N$ with $h_n^N \in \mathcal{C}^\infty_0(D_s)$ of the canonical form \eqref{eq:hxy}, and
    \item $P_n^N = (\mathrm{Id} + p_n^N) \circ L_n$ with $p_n^N \in \mathcal{R}_+^{(\ell)}$
\end{itemize}
for $0 \le n \le N-1$. We construct these maps so that the following diagram commutes up to $\ell$-jet equivalence in the directions transversal to $E_s$, with $H_N^N = \mathrm{Id}$:
\begin{equation}\label{cd:N}
\begin{CD}
U @>{f_0}>> U @>{f_1}>> \dots @>{f_{N-2}}>> U @>{f_{N-1}}>> U\\
@VV{H_0^N}V @VV{H_1^N}V @. @VV{H_{N-1}^N}V @VV{H_N^N= \mathrm{Id}}V\\
\real^d @>{P_0^N}>> \real^d @>{P_1^N}>> \dots @>{P_{N-2}^N}>> \real^d @>{P_{N-1}^N}>> \real^d
\end{CD}
\end{equation}
This is achieved by backward induction on $n$ for $0 \le n < N$. Suppose that $H_{n+1}^N$ has been defined. We determine $p_n^N$ and $h_n^N$ as follows:

\medskip
\noindent
\textbf{Step $\mathbf{1_n}$:} We choose $p_n^N \in \mathcal{R}^{(\ell)}_+$ such that
\begin{equation}\label{eq:step1}
\bar{h}_{n+1}^N := \left[(\mathrm{Id} + p_n^N)^{-1} \circ (\mathrm{Id} + h_{n+1}^N) \circ (\mathrm{Id} + g_n) - \mathrm{Id}\right]^{(\ell)} \in \mathcal{C}^{(\ell)}_0(L_n(D_s))
\end{equation}
belongs to the subspace $\overline{\mathcal{R}}^{(\ell)}_-(L_n(D_s))$.

\medskip
\noindent
\textbf{Step $\mathbf{2_n}$:} We then choose $h_n^N \in \mathcal{C}^{\infty}_0(D_s)$ of the form \eqref{eq:hxy} so that $[h_n^N]^{(\ell)} = \mathcal{L}_n(\bar{h}_{n+1}^N) \in \overline{\mathcal{R}}_-^{(\ell)}(D_s)$, where $\mathcal{L}_n$ is the conjugation operator associated with $L_n$.
\medskip

To ensure that this inductive procedure is well-defined, we verify that the choice in Step $1_n$ is possible at each step. Recall that we are assuming \eqref{eq:gn} for a small constant $\delta > 0$. If we also assume that $\|[h_{n+1}^N]^{(\ell)}\|_{C^R} < \delta'$ for some small $\delta' > 0$, then by the Implicit Function Theorem and the direct sum decomposition \eqref{eq:prod}, there exists a unique $p_n^N \in \mathcal{R}_+^{(\ell)}$ satisfying \eqref{eq:step1}. Thus, Step $1_n$ is valid as long as $\|[h_{n+1}^N]^{(\ell)}\|_{C^R} < \delta'$. Furthermore, by the normalization of the projection $\pi$ with respect to the $C^R$ norm, for any arbitrarily small constant $\eta > 0$, we may take $\delta$ and $\delta'$ so small that
\[
\|\bar{h}_{n+1}^N\|_{C^R} \le (1+\eta)(\delta + \delta').
\]
In Step $2_n$, the contraction property \eqref{eq:contraction} of $\mathcal{L}_n$ on $\overline{\mathcal{R}}^{(\ell)}_-(L_n(D_s))$ yields
\[
\|[h_n^N]^{(\ell)}\|_{C^R} \le \rho_R \|\bar{h}_{n+1}^N\|_{C^R} \le \rho_R (1+\eta) (\delta + \delta').
\]
Therefore, by choosing $\eta$, $\delta$, and $\delta'$ so that $(1+\eta)\rho_R < 1$ and
\[
\rho_R (1+\eta)(\delta + \delta') < \delta',
\]
the backward induction proceeds with the uniform estimate $\|[h_n^N]^{(\ell)}\|_{C^R} \le \delta'$ for all $0 \le n \le N$. It then follows that $\|p_n^N\|_{C^R}$ is also uniformly bounded.

Next, we fix $n \ge 0$ and examine the dependence of $h_n^{N}$ on $N \ge n$. Note that $h_n^{N+1}$ is obtained by replacing the terminal condition $H_N^N = \mathrm{Id}$ with $\mathrm{Id} + h_N^{N+1}$ at step $N$ of the backward inductive construction. Since each step of the induction contracts differences in the $C^R$ norm $\|\cdot\|_{C^R}$ by a factor of at most $(1+\eta)\rho_R < 1$, we obtain
\[
\|[h_n^{N+1} - h_n^{N}]^{(\ell)}\|_{C^R} \le ((1+\eta)\rho_R)^{N-n} \delta'.
\]
Consequently, as $N \to \infty$, $h_n^N$ converges in the $C^R$ norm to a limit $h_n$, and $p_n^N$ likewise converges to a limit $p_n \in \mathcal{R}_+^{(\ell)}$. Setting
\[
H_n = \mathrm{Id} + h_n = \lim_{N \to \infty} (\mathrm{Id} + h_n^N), \quad P_n = (\mathrm{Id} + p_n) \circ L_n = \lim_{N \to \infty} (\mathrm{Id} + p_n^N) \circ L_n,
\]
we obtain the following diagram, which commutes up to $\ell$-jet equivalence:
\begin{equation}\label{cd:N2}
\begin{CD}
    U @>{f_0}>> U @>{f_1}>> U @>{f_2}>> U @>{f_3}>> \dots \\
    @VV{H_0}V @VV{H_1}V @VV{H_2}V @VV{H_3}V \\
    \real^d @>{P_0}>> \real^d @>{P_1}>> \real^d @>{P_2}>> \real^d @>{P_3}>> \dots
\end{CD}
\end{equation}
From the contraction argument above, one also sees that the $C^R$ maps $H_n$ are unique up to $\ell$-jet equivalence and composition on the left with elements of $\mathcal{G}_+^{(\ell)}$ (we leave the straightforward details to the reader). Finally, since $R$ can be chosen arbitrarily large and the limit $[h_n]^{(\ell)}$ is independent of $R$, the maps $H_n$ are in fact of class $C^\infty$, which completes the proof.

\subsection{Proof of Theorem \ref{th:mainX}}
In what follows, we assume that the reader is familiar with the basic results of Pesin theory, particularly the construction of Lyapunov adapted metrics and charts (see \cite{Pesin, TsuR}). Let $\exp_p \colon T_pM \to M$ denote the exponential map at $p \in M$. For $n \ge 0$, we consider the local representation of the dynamics:
\[
\tilde{f}_n := \exp^{-1}_{f^{n+1}(p)} \circ f \circ \exp_{f^n(p)} \colon D \to T_{f^{n+1}(p)}M,
\]
where $D \subset T_{f^n(p)}M$ is an open neighborhood of the origin of uniform radius for all $n \ge 0$.

Next, we choose a sequence of linear isomorphisms $\psi_{n} \colon \mathbb{R}^d \to T_{f^n(p)}M$ for $n \ge 0$ such that:
\begin{itemize}
    \item $\psi_n(E_k) = Df_p^n(E_k(p))$ for $1 \le k \le m$;
    \item $\|\psi_n\| \le 1$, and $\|\psi_n^{-1}\|$ grows sub-exponentially as $n \to \infty$.
\end{itemize}
Since the point $p$ is forward Lyapunov regular, one can choose the sequence $\{\psi_n\}_{n \ge 0}$ so that the linearizations $L_n := D(\psi_{n+1}^{-1} \circ \tilde{f}_n \circ \psi_n)_0$ satisfy the spectral condition \eqref{eq:spec}.

Furthermore, by Pesin's stable manifold theorem, after rescaling closely around the origin (at a sub-exponential rate in $n$), we can modify the coordinate charts $\exp_{f^n(p)} \circ \psi_n$ by composing them with a sequence of $C^\infty$ diffeomorphisms so that the resulting local representations $f_n$ satisfy all the assumptions specified in Section~\ref{sec:setting}, including \eqref{ass:es} and the uniform boundedness of $\{f_n\}_{n \ge 0}$ in the $C^\infty$ topology. Thus, applying Theorem \ref{th:abs} to the sequence $\{f_n\}_{n \ge 0}$, we obtain the existence and uniqueness of the normal coordinates as claimed in Theorem \ref{th:mainX}.

\bibliographystyle{amsplain}
\bibliography{myBib.bib}
\end{document}

%% file: pict1.tex
\begin{tikzpicture}[scale=1.2,>=Stealth, line cap=round]

  \draw[thick] (-2,0) -- (2,0);     
  \draw[thick] (0,-1.5) -- (0,1.5); 

  \fill (0,0) circle (2pt);

  \draw[thick,->] (-2,0) -- (-1,0);
  \draw[thick,->] ( 2,0) -- ( 1,0);


  \node at (1.8, 0.35) {$\mathcal{R}_-$};
  \node at (0.95, 1.6) {${\mathcal{R}}_+$};
  \node at (-0.3,-0.3) {0};

\end{tikzpicture}

%% file: pict2.tex

\begin{tikzpicture}[
  line cap=round,
  x={(1.00cm,-0.35cm)},   
  y={(0.85cm, 0.18cm)},   
  z={(0.00cm, 1.00cm)}    
]

  \draw[thick,->] (0,0,0) -- (4.5,0,0) node[anchor=west] {$W^u_{loc}(p)$};
  \draw[thick,->] (0,0,0) -- (0,3.0,0) node[anchor=south] {$E^s$};
  \draw[thick,->] (0,0,0) -- (0,0,1.6) node[anchor=south] {$E^c$};

  \fill (4,0,0) circle (1.2pt);
\node[below left] at (4,0,0) {$q$};
\node[below left] at (0,0,0) {$p$};
\node[below right] at (4,3.2,0.2) {$W^s_{loc}(q)$};
  \fill (0,0,0) circle (1.2pt);


  \fill (0.400,0,0) circle (0.8pt);
  \draw (0.400,0,0)
    .. controls (0.6879447283,0.700,0.0989592277) and (0.8574370131,1.400,0.0214584583)
    .. (1.1409028595,2.200,0.2016094424);

  \fill (0.800,0,0) circle (0.8pt);
  \draw (0.800,0,0)
    .. controls (1.0878389271,0.700,0.0004084603) and (1.3449871725,1.400,0.0015365829)
    .. (1.6343153836,2.200,0.1762150972);

  \fill (1.200,0,0) circle (0.8pt);
  \draw (1.200,0,0)
    .. controls (1.4433052647,0.700,0.1003499421) and (1.7471590976,1.400,0.1043890840)
    .. (1.9844818738,2.200,0.1282174527);

  \fill (1.600,0,0) circle (0.8pt);
  \draw (1.600,0,0)
    .. controls (1.8839953550,0.700,0.2287137391) and (2.1648352555,1.400,0.0116891627)
    .. (2.4114395135,2.200,0.1825800158);

  \fill (2.000,0,0) circle (0.8pt);
  \draw (2.000,0,0)
    .. controls (2.2277445082,0.700,0.2519201916) and (2.5879436393,1.400,0.3519461036)
    .. (2.8284175300,2.200,0.3066805402);

  \fill (2.400,0,0) circle (0.8pt);
  \draw (2.400,0,0)
    .. controls (2.6558123767,0.700,0.2045809366) and (2.9193966178,1.400,0.3318875080)
    .. (3.2046508794,2.200,0.0233569992);

  \fill (2.800,0,0) circle (0.8pt);
  \draw (2.800,0,0)
    .. controls (3.0905886140,0.700,0.2311197927) and (3.2538200108,1.400,0.2319041585)
    .. (3.5975604979,2.200,0.2793631794);

  \fill (3.200,0,0) circle (0.8pt);
  \draw (3.200,0,0)
    .. controls (3.4827195259,0.700,0.0105343003) and (3.6952275903,1.400,0.4741802337)
    .. (3.9632862276,2.200,0.2157161904);

  \fill (3.600,0,0) circle (0.8pt);
  \draw (3.600,0,0)
    .. controls (3.8775886634,0.700,0.1233160418) and (4.0893148927,1.400,0.2899684508)
    .. (4.3570243802,2.200,0.1948460482);

  \fill (4.000,0,0) circle (0.8pt);
  \draw (4.000,0,0)
    .. controls (4.2656484314,0.700,0.1529677666) and (4.5797274560,1.400,0.5110650360)
    .. (4.7521805210,2.200,0.2331915917);

\end{tikzpicture}

%% file: myBib.bib
@misc{Brown,
    Author = {Aaron Brown and Alex Eskin and Simion Filip and Federico Rodriguez Hertz},
    Title = {Normal forms for contracting dynamics, revisited},
    Year = {2024},
    archivePrefix = {arXiv},
    Note = {arXiv:2405.16208},
    Eprint = {arXiv:2405.16208},
}

@incollection{Kalinin,
    AUTHOR = {Kalinin, Boris},
     TITLE = {Non-stationary normal forms for contracting extensions},
 BOOKTITLE = {A vision for dynamics in the 21st century---the legacy of
              {A}natole {K}atok},
     PAGES = {207--231},
 PUBLISHER = {Cambridge Univ. Press, Cambridge},
      YEAR = {2024},
      ISBN = {978-1-009-27890-4},
   MRCLASS = {37C86 (37C15 37D99 37G99)},
  MRNUMBER = {4685044},
}

@article{KS1,
    AUTHOR = {Kalinin, Boris and Sadovskaya, Victoria},
     TITLE = {Normal forms for non-uniform contractions},
   JOURNAL = {J. Mod. Dyn.},
  FJOURNAL = {Journal of Modern Dynamics},
    VOLUME = {11},
      YEAR = {2017},
     PAGES = {341--368},
      ISSN = {1930-5311,1930-532X},
   MRCLASS = {37D10 (34C20 37D20 37D25)},
  MRNUMBER = {3642250},
MRREVIEWER = {Viorel\ Ni\c tic\u a},
       DOI = {10.3934/jmd.2017014},
       URL = {https://doi.org/10.3934/jmd.2017014},
}

@article{K,
    AUTHOR = {Kalinin, Boris and Sadovskaya, Victoria},
     TITLE = {Normal forms on contracting foliations: smoothness and
              homogeneous structure},
   JOURNAL = {Geom. Dedicata},
  FJOURNAL = {Geometriae Dedicata},
    VOLUME = {183},
      YEAR = {2016},
     PAGES = {181--194},
      ISSN = {0046-5755,1572-9168},
   MRCLASS = {37C15 (34C20 37D10 37D30)},
  MRNUMBER = {3523125},
MRREVIEWER = {Shigenori\ Matsumoto},
       DOI = {10.1007/s10711-016-0153-5},
       URL = {https://doi.org/10.1007/s10711-016-0153-5},
}

@article{Melnick,
    AUTHOR = {Melnick, Karin},
     TITLE = {Non-stationary smooth geometric structures for contracting
              measurable cocycles},
   JOURNAL = {Ergodic Theory Dynam. Systems},
  FJOURNAL = {Ergodic Theory and Dynamical Systems},
    VOLUME = {39},
      YEAR = {2019},
    NUMBER = {2},
     PAGES = {392--424},
      ISSN = {0143-3857,1469-4417},
   MRCLASS = {37C60 (37C40 37G05 37H15)},
  MRNUMBER = {3893265},
MRREVIEWER = {Mahesh\ Nerurkar},
       DOI = {10.1017/etds.2017.38},
       URL = {https://doi.org/10.1017/etds.2017.38},
}

@article{Oseledec,
    AUTHOR = {Oseledec, V. I.},
     TITLE = {A multiplicative ergodic theorem. {C}haracteristic {L}japunov,
              exponents of dynamical systems},
   JOURNAL = {Trudy Moskov. Mat. Ob\v s\v c.},
  FJOURNAL = {Trudy Moskovskogo Matemati\v ceskogo Ob\v s\v cestva},
    VOLUME = {19},
      YEAR = {1968},
     PAGES = {179--210},
      ISSN = {0134-8663},
   MRCLASS = {28.70 (34.00)},
  MRNUMBER = {240280},
MRREVIEWER = {J\'ozsef\ Sz\H ucs},
}

@misc{G,
    Author = {Andrey Gogolev and Martin Leguil and Federico Rodriguez Hertz},
    Title = {Smooth rigidity for 3-dimensional dissipative {A}nosov flows},
    Year = {2025},
    Note = {arXiv:2510.23872},
    Eprint = {arXiv:2510.23872},
}

@article{Pesin,
    AUTHOR = {Pesin, Ja.\ B.},
     TITLE = {Characteristic {L}japunov exponents, and smooth ergodic
              theory},
   JOURNAL = {Uspehi Mat. Nauk},
  FJOURNAL = {Akademija Nauk SSSR i Moskovskoe Matemati\v ceskoe Ob\v s\v
              cestvo. Uspehi Matemati\v ceskih Nauk},
    VOLUME = {32},
      YEAR = {1977},
    NUMBER = {4(196)},
     PAGES = {55--112, 287},
      ISSN = {0042-1316},
   MRCLASS = {34D05 (28A65 58F10 58F15)},
  MRNUMBER = {466791},
MRREVIEWER = {A.\ Morimoto},
}

@article{T,
    AUTHOR = {Tsujii, Masato},
     TITLE = {Exponential mixing for generic volume-preserving {A}nosov
              flows in dimension three},
   JOURNAL = {J. Math. Soc. Japan},
  FJOURNAL = {Journal of the Mathematical Society of Japan},
    VOLUME = {70},
      YEAR = {2018},
    NUMBER = {2},
     PAGES = {757--821},
      ISSN = {0025-5645,1881-1167},
   MRCLASS = {37A25 (37D20)},
  MRNUMBER = {3787739},
MRREVIEWER = {Saadet\ \"Oyk\"u\ Yurtta\c s},
       DOI = {10.2969/jmsj/07027595},
       URL = {https://doi.org/10.2969/jmsj/07027595},
}

@article{TZ,
    AUTHOR = {Tsujii, Masato and Zhang, Zhiyuan},
     TITLE = {Smooth mixing {A}nosov flows in dimension three are
              exponentially mixing},
   JOURNAL = {Ann. of Math. (2)},
  FJOURNAL = {Annals of Mathematics. Second Series},
    VOLUME = {197},
      YEAR = {2023},
    NUMBER = {1},
     PAGES = {65--158},
      ISSN = {0003-486X,1939-8980},
   MRCLASS = {37D20 (37A25 37C30 53D25)},
  MRNUMBER = {4513143},
MRREVIEWER = {Ana\ Rechtman},
       DOI = {10.4007/annals.2023.197.1.2},
       URL = {https://doi.org/10.4007/annals.2023.197.1.2},
}

@article{GK,
    AUTHOR = {Guysinsky, M. and Katok, A.},
     TITLE = {Normal forms and invariant geometric structures for dynamical
              systems with invariant contracting foliations},
   JOURNAL = {Math. Res. Lett.},
  FJOURNAL = {Mathematical Research Letters},
    VOLUME = {5},
      YEAR = {1998},
    NUMBER = {1-2},
     PAGES = {149--163},
      ISSN = {1073-2780},
   MRCLASS = {58F36 (58F15)},
  MRNUMBER = {1618331},
MRREVIEWER = {Nantian\ Qian},
       DOI = {10.4310/MRL.1998.v5.n2.a2},
       URL = {https://doi.org/10.4310/MRL.1998.v5.n2.a2},
}

@misc{Filip,
    Author = {Simion Filip},
    Title = {Measure Rigidity beyond Homogeneous Dynamics},
    Year = {2025},
    Eprint = {arXiv:2512.13865},
    Note = {arXiv:2512.13865},
}

@misc{GLeclerc,
    Author = {Ga\"etan Leclerc},
    Title = {Fourier decay of equilibrium states and the {F}ibonacci {H}amiltonian},
    Year = {2025},
    Eprint = {arXiv:2507.23731},
    Note = {arXiv:2507.23731},
}

@article{TsuR,
    AUTHOR = {Tsujii, Masato},
     TITLE = {Regular points for ergodic {S}ina\u i\ measures},
   JOURNAL = {Trans. Amer. Math. Soc.},
  FJOURNAL = {Transactions of the American Mathematical Society},
    VOLUME = {328},
      YEAR = {1991},
    NUMBER = {2},
     PAGES = {747--766},
      ISSN = {0002-9947,1088-6850},
   MRCLASS = {58F11},
  MRNUMBER = {1072103},
MRREVIEWER = {A.\ I.\ Danilenko},
       DOI = {10.2307/2001802},
       URL = {https://doi.org/10.2307/2001802},
}

@misc{EPZ,
    Author = {Alex Eskin and Rafael Potrie and Zhiyuan Zhang},
    Title = {Geometric properties of partially hyperbolic measures and applications to measure rigidity},
    Year = {2023},
    Eprint = {arXiv:2302.12981},
    Note = {arXiv:2302.12981},
}

@article{Guy,
    AUTHOR = {Guysinsky, M.},
     TITLE = {The theory of non-stationary normal forms},
   JOURNAL = {Ergodic Theory Dynam. Systems},
  FJOURNAL = {Ergodic Theory and Dynamical Systems},
    VOLUME = {22},
      YEAR = {2002},
    NUMBER = {3},
     PAGES = {845--862},
      ISSN = {0143-3857,1469-4417},
   MRCLASS = {37C15 (37G05)},
  MRNUMBER = {1908557},
MRREVIEWER = {Boris\ Hasselblatt},
       DOI = {10.1017/S0143385702000421},
       URL = {https://doi.org/10.1017/S0143385702000421},
}
